\documentclass[11pt]{article}
\usepackage{latexsym, amsfonts, tikz}  
\newcommand{\be}{\begin{equation}}
\newcommand{\ee}{\end{equation}}
\newcommand{\bea}{\begin{eqnarray}}
\newcommand{\eea}{\end{eqnarray}}
\newcommand{\bean}{\begin{eqnarray*}}
\newcommand{\eean}{\end{eqnarray*}}
\newcommand{\brray}{\begin{array}}
\newcommand{\erray}{\end{array}}
\newcommand{\ben}{\begin{equation}{nonumber}}
\newcommand{\een}{\end{equation}{nonumber}}

\newtheorem{dfn}{Definition}[section]
\newtheorem{thm}[dfn]{Theorem}
\newtheorem{lmma}[dfn]{Lemma}
\newtheorem{ppsn}[dfn]{Proposition}
\newtheorem{crlre}[dfn]{Corollary}
\newtheorem{xmpl}[dfn]{Example}
\newtheorem{rmrk}[dfn]{Remark}
\newcommand{\bdfn}{\begin{dfn}}
\newcommand{\bthm}{\begin{thm}}
\newcommand{\blmma}{\begin{lmma}}
\newcommand{\bppsn}{\begin{ppsn}}
\newcommand{\bcrlre}{\begin{crlre}}
\newcommand{\bxmpl}{\begin{xmpl}}
\newcommand{\brmrk}{\begin{rmrk}}
\newcommand{\edfn}{\end{dfn}}
\newcommand{\ethm}{\end{thm}}
\newcommand{\elmma}{\end{lmma}}
\newcommand{\eppsn}{\end{ppsn}}
\newcommand{\ecrlre}{\end{crlre}}
\newcommand{\exmpl}{\end{xmpl}}
\newcommand{\ermrk}{\end{rmrk}}

\newcommand{\IC}{\mathbb{C}}

\newcommand{\IR}{\mathbb{R}}

\newcommand{\IZ}{\mathbb{Z}}

\newcommand{\cla}{{\cal A}}
\newcommand{\clb}{{\cal B}}
\newcommand{\clc}{{\cal C}}
\newcommand{\cld}{{\cal D}}
\newcommand{\cle}{{\cal E}}
\newcommand{\clf}{{\cal F}}
\newcommand{\clg}{{\cal G}}
\newcommand{\clh}{{\cal H}}
\newcommand{\cli}{{\cal I}}

\newcommand{\clk}{{\cal K}}
\newcommand{\cll}{{\cal L}}
\newcommand{\clm}{{\cal M}}

\newcommand{\clp}{{\cal P}}
\newcommand{\clq}{{\cal Q}}

\newcommand{\cls}{{\cal S}}

\newcommand{\clu}{{\cal U}}
\newcommand{\clv}{{\cal V}}
\newcommand{\clw}{{\cal W}}

\def\a*{{\cal A}_{h,*}}
\def\B{{\cal B}(h)}
\def\B1{{\cal B}_1(h)}
\def\b{{\cal B}^{\rm s.a.}(h)}
\def\b1{{\cal B}^{\rm s.a.}_1(h)}

\newcommand{\ot}{\otimes}

\newcommand{\raro}{\rightarrow}

\def \qed {$\Box$}

\begin{document}
\begin{center}
{\large {\bf Existence and examples of quantum isometry groups  for a class of compact metric spaces}}\\
by\\ 
{\large Debashish Goswami }\footnote {Work partially supported by IUSSTF (through Indo-US fellowship),  
Dept. of Science and Technology  of Govt. of India  (Swarnajayanti fellowship and project) and  the Indian National Science Academy.}\\ 
{\large Stat-Math Unit, Kolkata Centre,}\\
{\large Indian Statistical Institute}\\
{\large 203, B. T. Road, Kolkata 700 108, India}\\
{e mail: goswamid@isical.ac.in}\\
\end{center}

\begin{abstract}
 We  formulate a definition of isometric action of a compact quantum group (CQG) on a compact metric space, 
  generalizing   Banica's definition for finite metric spaces. For metric spaces $(X,d)$ which can be isometrically embedded in some Euclidean space, we 
   prove the existence of a universal object in the category of the 
  compact quantum groups  acting isometrically on $(X,d)$.  In fact, our existence theorem applies to a larger class, namely for any compact metric 
   space $(X,d)$ which admits a one-to-one continuous map $f : X \raro \IR^n$ for some $n$ such that $d_0(f(x),f(y))=\phi(d(x,y))$ (where $d_0$ is the Euclidean metric) 
   for some homeomorphism $\phi$ of $\IR^+$. 
  
  As concrete examples, 
   we obtain Wang's quantum permutation group $\cls_n^+$ and also the free wreath product of $\IZ_2$ by $\cls_n^+$ as
    the quantum isometry groups for certain compact  connected metric spaces constructed  by taking topological joins of intervals in \cite{huang1}. 
\end{abstract}
\section{Introduction}
 It is a very natural and interesting question to study quantum symmetries of classical spaces, particularly metric spaces. 
In fact, motivated by some suggestions of Alain Connes, S. Wang defined (and proved existence) of compact quantum group analogues of the classical 
symmetry or automorphism groups of various types of finite structures such as finite sets  
 and finite dimensional matrix algebras (see \cite{wang1}, \cite{wang2}), 
and then these quantum groups were investigated in depth by a number of mathematicians including Wang, Banica, Bichon and others 
(see, for example, \cite{ban1},  \cite{bichon} and the references therein). We should also mention \cite{manin_book} for an algebraic 
pre-cursor of Wang's universal compact quantum groups. One may expect that in general, there will be many more 
 quantum symmetries of a given classical space than its classical group
     symmetries.  Indeed, this is the case for a finite set of size $n$ with  
      $n \geq 4$, which  has an infinite dimensional 
      compact quantum group (`quantum permutation group' $\cls_n^+$ a la Wang) of symmetries.

However, it is important to extend these 
ideas and constructions to the `continuous' or `geometric' set-up. In  a series of articles initiated by us in \cite{goswami} and followed up
 by several other mathematicians, e.g. Bhowmick,  Banica, Soltan,  De-Commer, Thibault, Das, Joardar, Skalski, Mandal, just to name a few, 
 we have formulated, studied and computed  quantum group analogues 
of the group of isometries  of Riemannian manifolds, including in fact noncommutative geometric set-up in the sense 
of \cite{connes} as well. 
  It is a natural question  whether such construction can be done in a purely metric space set-up, without assuming any finer geometric 
(e.g. Riemannian or spin) structures. 
 
 In fact, there is another strong reason for considering metric space set-up, coming from our computations of quantum isometry groups of Riemannian manifolds as well as 
  the general results obtained in \cite{nonexistence}, which we now explain. Although it is easy to construct  faithful isometric actions of 
   genuine (noncommutative as a $C^*$ algebra) compact quantum groups on disconnected compact Riemannian manifolds having at least four identical components, using 
    the  natural action of  quantum permutation group $\cls^+_4$, the situation dramatically changes if we look for such actions on compact connected smooth manifolds. 
     Indeed, it is proved in \cite{nonexistence} that there cannot be any faithful smooth action of a genuine compact quantum group on 
      a compact connected smooth manifold, which in particular implies that the quantum isometry group for any such  manifold M must coincide with $C(ISO(M))$. 
     However, Huang (\cite{huang1})  constructed several examples of faithful actions of genuine compact quantum groups 
      on compact connected spaces which are obtained by topologically gluing smooth manifolds, e.g. topological join of compact intervals. In \cite{etingof}, 
       an example of finite-dimensional genuine compact quantum group acting faithfully on the space of a non-smooth algebraic variety has been given. Thus, there seem to 
        be many interesting actions of genuine compact quantum groups on connected manifolds with singularities, typically embedded in  Euclidean spaces. The actions 
         in the above-mentioned examples are also isometric in a natural sense if we consider the restriction of Euclidean metric. A general framework for studying 
       such actions should be obtained if we can formulate a satisfactory theory of quantum isometry groups for metric spaces, covering at least those which are embedded isometrically 
       in          Euclidean spaces.

This aim  is achieved in the present article. It may be mentioned that   we have proposed in \cite{qiso_disc} a natural definition of `isometric' action of a (compact) 
quantum group on an arbitrary 
compact metric space, extending Banica's definition of quantum isometry group of finite metric spaces, and showed in some explicit examples the existence 
of a universal object in the category of all such compact quantum groups acting isometrically on a given metric space. 
There is also an attempt in \cite{sabbe} to give such a formulation for even more general framework of compact quantum metric spaces a la Rieffel (\cite{rieffel}).
However, the formulation of \cite{sabbe} does not seem to be very convenient for computations as it involves inequalities rather than equalities in the definition of `isometric action'.
In the present paper we slightly modify the definition proposed by us in \cite{qiso_disc} (for finite spaces it is still the same) and then 
prove the existence of a universal object in the category of all compact quantum groups acting isometrically on a given space, for a quite large   class of 
metric  spaces, including (but not limited to)  those which are isometrically embeddable in some $\IR^n$. In fact, our existence theorem applies to any compact metric 
   space $(X,d)$ which admits a one-to-one continuous map $f : X \raro \IR^n$ for some $n$ such that $d_0(f(x),f(y))=\phi(d(x,y))$ 
   (where $d_0$ is the Euclidean metric) for some homeomorphism $\phi$ of $\IR^+$. This class also includes all finite metric spaces, so we do extend Banica's 
   theory of quantum isometry groups.  Let us also mention the recent paper \cite{alex_preprint} by Chirvasitu for comparison of  various definitions of isometric actions.
    
We begin with a brief account of compact quantum groups and their actions in Section 2. Then we formulate the definition of isometric actions by compact quantum groups 
 on compact metric spaces in Section 3, followed by the main results about  existence of the quantum isometry groups in Section 4. 
 We conclude the article with some  explicit computations of  quantum isometry groups for a number of interesting concrete examples in Section 5.

\section{Quantum groups and their actions}
Let us very briefly review the basics of compact quantum groups, their actions and representations, referring the reader to  \cite{woro1}, \cite{woro2}, \cite{podles}  for details (see also 
 \cite{alex_preprint}, \cite{huang2}).
A compact quantum group (CQG for short) is a  unital $C^*$ algebra $\cls$ with a unital $\ast$-homomorphism 
 $\Delta$ from $\cls$ to $\cls \ot \cls$ (injective tensor product) 
which is coassociative, i.e. $(\Delta \ot {\rm id})\circ \Delta=({\rm id} \ot \Delta) \circ \Delta$ and 
each of the linear spans of $\Delta(\cls)(\cls \ot 1)$ and that of $\Delta(\cls) (1 \ot \cls)$ is norm-dense in $\cls \ot \cls$. 
From this condition, one can obtain a canonical dense unital $\ast$-subalgebra $\cls_0$ of $\cls$ and linear maps $\kappa$ and 
$\epsilon$ (called the antipode and the counit respectively) defined on $\cls_0$ making it a Hopf $\ast$-algebra. 
 Any CQG $\cls$ has a unique   state $h$ called the Haar state satisfying $({\rm id} \ot h)\circ \Delta(a)=(h \ot {\rm id}) \circ \Delta(a)=h(a)1$ 
 for all $a$. The Haar state need not be faithful in general, though it is always faithful on $\cls_0$ at least. 
 The CQG $\cls_r:=\pi_r(\cls)$, where $\pi_r: \cls \raro \clb(L^2(h))$ is the GNS representation, is called the reduced CQG corresponding to $\cls$. 

 For $T \in \cll(\clh \ot \cls)$ (i.e. adjointable $\cls$-linear map)  we write $T_{23}:=I_\clh \ot T$, $T_{13}:=\sigma_{12} \circ T_{23} \circ \sigma_{12}$,
  $T^{12}:=T \ot {\rm id}_\cls$ and $T^{13}:=\sigma_{23} \circ T^{12} \circ \sigma_{23}$, where 
   $\sigma_{ij}$ denotes the map which flips $i$-th and $j$-th copies of tensor components. Now, a unitary representation of   $(\cls, \Delta)$ in a Hilbert space $\clh$ is  a complex linear map  $U$ from  $\clh$ to 
$\clh \ot \cls$, such that the $\cls$-linear map $\tilde{U}$ on $\clh \ot \cls$ defined by $\tilde{U}(\xi \ot b)=U(\xi) b$, 
 ( $\xi \in \clh, b \in \cls$) extends to a well-defined $\cls$-linear unitary  satisfying $({\rm id} \ot \Delta) (\tilde{U})=\tilde{U}^{12}\tilde{U}^{13}$. 
A CQG is called a compact matrix quantum group if  it has a dense $\ast$-subalgebra generated by the matrix coefficients of a 
 finite dimensional unitary representation (called the fundamental representation). Most (but not all) CQG's considered by us will be of this type.

A (co)action of a CQG $(\cls, \Delta)$ on  a unital $C^*$ algebra $\clc$ is a unital $C^*$-homomorphism 
$\beta : \clc \raro \clc \ot \cls$ such that $(\beta \ot {\rm id}) \circ \beta=({\rm id} \ot \Delta) \circ \beta$ and the linear span of 
$\{ \beta(\clc)(1 \ot \cls)\}$ is norm-dense in $\clc \ot \cls$. Given such  an action $\beta$, one can find  a unital dense $\ast$-subalgebra $\clc_0$ of $\clc$, 
called the spectral subalgebra,  such that $\beta$ maps $\clc_0$ into $\clc_0 \ot_{\rm alg} \cls_0$ 
  and we also have 
   $({\rm id} \ot \epsilon)\circ
 \beta={\rm id}$ on $\clc_0$. We'll occasionally use the standard Sweedler notation for coproduct and action of Hopf algebra, i.e. write $\Delta(q)=q_{(1)} \ot q_{(2)},$ 
  $\beta(a)=a_{(0)} \ot a_{(1)}.$

  
  
 %
 We say that the action $\beta$ is faithful if the $\ast$-subalgebra of $\cls$ generated by elements of the form $(\omega \ot {\rm id})(\beta(a))$, 
 where $a \in \clc$, $\omega$ being a bounded linear functional on $\clc$, is norm-dense in $\cls$. We refer to the Proposition 3.2 of \cite{huang2} for various equivalent
  descriptions of faithfulness. 
 %
  We also recall (Theorem 3.23, \cite{huang2}) the well-known fact that if a CQG $\cls$ which acts faithfully on $C(X)$, where $X$ is a compact Hausdorff space, 
   then $\cls$ is of Kac type, i.e. the Haar state is tracial and the antipode $\kappa$ admits a norm-bounded extension on $\cls_r$ satisfying $\kappa^2={\rm id}.$ For any such 
    action $\beta$ on $C(X)$, we call a probability measure $\mu$ on $X$ $\beta$-invariant if the corresponding state $\phi$ (say) is $\beta$-invariant, 
    i.e. $(\phi \ot {\rm id})(\beta(f))=\phi(f) 1\equiv (\int_X fd\mu)1$ $\forall f$.

Let us give  two examples of CQG (in fact compact matrix quantum groups). The first  is the quantum permutation group $\cls_n^+$ due to  Wang, 
 which is the universal unital $C^*$ algebra generated by $n^2$ orthogonal projections $\{ w_{ij} \}$, satisfying $ \sum_{i=1}^n w_{ij}=1,~~\sum_{j=1}^n w_{ij}=1.$
  The coproduct is defined by $\Delta(w_{ij})=\sum_k w_{ik} \ot w_{kj}.$  
  
  The other example is 
   the quantum free orthogonal group $A_o(n)$, which is defined to be the universal unital $C^*$ algebra generated by self-adjoint elements $\{ q^o_{ij} \}$ subject to the conditions 
    that $(( q^o_{ij} ))$ as well as $(( q^o_{ji} ))$ are unitaries.  The coproduct is given by $\Delta(q^o_{ij})=\sum_k q^o_{ik} \ot q^o_{kj}.$

\section{Definition of isometric action of compact quantum groups}

\bdfn 
\label{isodef}
Given an  action $\beta$ of a CQG $\cls$ on  $\clc=C(X)$ (where $(X,d)$ is a compact metric space), we say that $\beta$ is  `isometric'  if the corresponding reduced action $\beta_r
:=({\rm id} \ot \pi_r)\circ \beta$ of $\cls_r$ satisfies the following:
\be \label{isodef_formula} ({\rm id}_\clc \ot \beta_r)(d)=\sigma_{23} \circ (({\rm id}_\clc \ot \kappa) \circ \beta_r \ot {\rm id}_\clc)(d),\ee
 where as before, $\sigma_{23}$ denotes the flip of the second and third tensor copies, $d \in C(X \times X)$ denotes the metric and $\kappa$ denotes the (bounded) antipode.
\edfn

 

  
\blmma
\label{iso_cond_lemma}
Given a  $C^*$-action $\beta$ of a CQG $\cls$ on $C(X)$, the following are equivalent:\\
(i) The action is isometric.\\
(ii) $\forall x, y \in X$, one has \be \label{iso_cond} \beta_r(d_x)(y)=\kappa(\beta_r(d_y)(x)),\ee where $d_x(z):=d(x,z)$.
\elmma
{\it Proof:}\\
The equivalence  of (i) and (ii) is a consequence of  the continuity of the map $x \mapsto d_x \in C(X)$, and hence
 (by the norm-contractivity of $\beta_r$), the continuity of $x \mapsto \beta_r(d_x) \in C(X) \ot \cls$.\qed\\
 
  The condition (\ref{iso_cond}) is termed as `D-isometry' in \cite{alex_preprint}. Moreover, it is proved in   \cite{alex_preprint}, Proposition 3.8 that
   for any isometric action $\beta$, $\beta_r$ is injective.
    Hence   by the discussion in Subsection 3.2 of  \cite{huang2}
      and Lemma 4.9 of \cite{nonexistence}, 
     there is a faithful (i.e. having full support) $\beta_r$-invariant Borel probability measure $\mu$ on $X$ 
      so that  $\beta_r$ extends to a unitary representation $U=U_\mu$ on $L^2(X,\mu)$, given by $U(f)=\beta_r(f)$. Viewing $d$ as a vector in $L^2(X \times X, \mu \times \mu),$
       it is easy to see that the condition (\ref{isodef_formula}) is equivalent to 
       $\tilde{U}_{23}(d \ot 1)=\tilde{U}_{13}^{-1}(d \ot 1)$. 
       Using the symmetry of $d$, i.e. $d=\sigma \circ d$, where $\sigma$ is the flip map  on $C(X) \ot C(X)$, this is also equivalent to 
        $\tilde{U}_{13}(d \ot 1)=\tilde{U}_{23}^{-1}(d \ot 1)$. Thus, an action is isometric if and only if  for some (hence any) invariant faithful $\mu$ as above, the 
         corresponding unitary $U$ satisfies
   \be \label{111} \tilde{U}_{13} \tilde{U}_{23}(d \ot 1)=d \ot 1=\tilde{U}_{23} \tilde{U}_{13}(d \ot 1). \ee
      Let $\tilde{W}:=\tilde{U}_{13}\tilde{U}_{23}$, $\tilde{Z}:=\tilde{U}_{23} \tilde{U}_{13}$ and denote by $\pi\equiv \pi_\mu : C(X) \raro 
   \clb(L^2(X,\mu))$ the GNS representation and let $\pi_\mu^{(2)}=\pi_\mu \ot \pi_\mu : C(X \times X) \raro \clb(L^2(X\times X,\mu \times \mu)).$  
   By straightforward calculations using the commutativity of $C(X)$ one can prove 
    $\tilde{U}^{-1}(\pi_\mu(f) \ot 1) \tilde{U}=(\pi_\mu \ot \kappa)(\beta_r(f)),$ for all $f$ belonging to the spectral subalgebra $C(X)_0$ and hence for all $f \in C(X)$.
     Then,  (\ref{isodef_formula}) is clearly equivalent to the following:
     $$ \tilde{W}(\pi_\mu^{(2)}(d) \ot 1)\tilde{W}^{-1}=\pi_\mu^{(2)}(d) \ot 1 =\tilde{Z}(\pi_\mu^{(2)}(d) \ot 1)\tilde{Z}^{-1}.$$



\brmrk
In case $\cls=C(G)$ and $\beta$ corresponds to a topological action of a compact group $G$ on $X$, it is clear that the definition of 
 isometric action  is nothing but the requirement $d(x, gy)=d(y,g^{-1}x) (=d(g^{-1}x,y)) $, which coincides with the usual definition of isometric group action.
Moreover, it follows from the equivalent condition (\ref{111}) that 
for a finite metric space $(X,d)$, the present definition does coincide with Banica's definition in \cite{ban1} as well as the one proposed in \cite{qiso_disc}. 

However,  it is not yet clear whether the  definition of isometric CQG-action in terms of an inequality given in \cite{sabbe} 
 is the same as the one given by us for a general compact metric space.
 For finite spaces the equivalence of the two definitions has been proved by them in \cite{sabbe} and it is proved in \cite{alex_preprint} that in general, 
 the definition given in this paper 
  implies that of  \cite{sabbe}. 
\ermrk

\section{Existence of a universal isometric action}
It is a natural question to ask: does there exist a universal object in the category (say ${\bf Q}_{X,d}$) of all CQG   
acting isometrically (in our sense) on $(X,d)$? For finite metric spaces, 
the answer is clearly affirmative, and the universal object is the quantum isometry group defined by Banica.
 We are not yet able to settle this question in full generality. However, we'll give an affirmative answer for a large class of metric spaces which can be isometrically embedded in 
 some $\IR^n$. The key idea is to prove that any isometric action for such metric spaces must be affine. One can compare this with 
  the classical situation: any group acting isometrically on a compact subset of $\IR^n$ (w.r.t. the Euclidean metric) must be a subgroup of 
   $O(n)$. In our more general framework, we have been able to obtain a similar result, with the free quantum orthogonal group $A_o(n)$ replacing the classical 
    orthogonal group $O(n)$. That is, any CQG acting isometrically on a compact subset of $\IR^n$ is shown to be a quantum subgroup of $A_o(n)$ in a canonical way. This allows us to prove 
     existence of a universal object as a suitable quotient of $A_o(n)$. 
  \subsection{Sufficient conditions for an action to be affine}
  The aim of this subsection is to prove 
   some  preparatory lemmas for showing (in the next subsection) that any isometric CQG action on a subset of an Euclidean space (with the induced metric) must be affine.

\blmma
 \label{iso_lin}
 Let $\clc$ be a unital $C^*$ algebra. Define $<< \cdot, \cdot >>$ on $\IR^m \ot_{\rm alg} \clc^{s.a.}$ by $<<Z,W>>:=\frac{1}{2}\sum_i (Z_iW_i+W_iZ_i)$. Let $F$ be a  
 function from $\IR^n$ to $\IR^m \ot_{\rm alg} \clc^{s.a.}$ which satisfies $<<F(x),F(y)>>=<x,y>1$ for all $x,y$, where $<\cdot,\cdot>$ is the Euclidean inner product of $\IR^n$. 
 Then $F$ must be linear.
 \elmma
 {\it Proof:}\\ 
 Let $\|A \|_\clc^2:=<<A,A>>$. It is easy to see that $\| A \|_\clc=0$ if and only if $A=0$ for $A \in \IR^m \ot_{\rm alg} \clc^{s.a.}$.  We now observe that $\| F(x+y)-F(x)-F(y)\|_\clc^2=0$ and $\| F(cx)-cF(x)\|_\clc^2=0$ by direct computation using the condition $<<F(x),F(y)>>=<x,y>1$.
 \qed\\
 \brmrk
 The bilinear form $<< \cdot, \cdot >>$ is not an $\clc^{s.a.}$ valued inner product in the sense of Hilbert module. It is only bilinear w.r.t. scalars, but not w.r.t. $\clc^{s.a.}$.
 \ermrk
 \blmma
 \label{iso_aff}
 Let $\clc$, $<<\cdot, \cdot>>$ and $\| \cdot \|_\clc^2$ be as in the statement and proof of the previous lemma, but let $F$ be a function from a nonempty subset $X$ of $\IR^m$ to
 $\IR^n \ot_{\rm alg} \clc^{s.a.}$ which satisfies $$ \| F(x)-F(y) \|_\clc^2=\| x-y\|^2 1$$ for all $x,y \in X$. Then $F$ is affine in the sense that there are elements
  $a_{ij}, i=1,...,n;j=1,...,m$ of $\clc^{s.a.}$ and $\xi=(\xi_1,...,\xi_n) \in
 \IR^n \ot_{\rm alg}  \clc^{s.a.}$ such that $$ F(x)=Ax+\xi,~~~A:=(( a_{ji})).$$
 \elmma
 {\it Proof:}\\
 Consider $x_0 $ in $X$ and the function $G:Y:= X-x_0 \equiv \{ x-x_0:~x \in X\} \raro \IR^n \ot_{\rm alg} \clc^{s.a.}$ by $G(y)=F(y+x_0)-F(x_0)$ for $y \in Y$. Then $\| G(y)-G(y^\prime)\|_\clc^2=\|y-y^\prime \|^21$ by hypothesis, and also $\| G(y)\|_\clc^2=\|y\|^21$. Observe that for any $Z,W \in \IR^n \ot_{\rm alg} \clc^{s.a.}$ we have $\| Z+W\|_\clc^2=\| Z\|_\clc^2+
 \|W\|_\clc^2+2<<Z,W>>$. Using this with $Z=G(y)$,$ W=G(y^\prime)$, we get that  \be \label{290}<<G(y),G(y^\prime)>>=<y,y^\prime>1~~~\forall y,y^\prime \in Y.\ee 
 This allows us to extend $G$ linearly to the span of $Y$ in $\IR^m$. Indeed, the extension  $G(\sum_i c_i y_i)=\sum_i c_i G(y_i)$, for $c_i \in \IR,$ $y_i \in Y$ 
 is well-defined because $\| \cdot \|_\clc^2$ of the right hand side equals the Euclidean norm square of $\sum_i c_i y_i$ by (\ref{290}). 
 We can then extend $G$ further on the whole of $\IR^m$ as a linear map denoted again by $G$. Thus, we get $A=(( a_{ji})),$ say, such that $G(y)=Ay$ for all $y \in \IR^m$. 
 This implies, $F(x)=G(x-x_0)+F(x_0)=Ax+\xi$, where $\xi=F(x_0)-Ax_0$. \qed\\
 We also get a slightly different criterion for the action to be affine.
\blmma
\label{affine_action}
Let $0 \in X \subseteq \IR^n$ and assume that the restriction of the coordinate functions of $\IR^n$ to $X$, say $X_1, \ldots, X_n$, are linearly independent. Let $\clc$ be a unital 
 $C^*$ algebra and $F_1,\ldots, F_n, G_1, \ldots, G_n$ be functions from $X$ to $\clc^{\rm s.a.}$ such that $$ \sum_i (F_i(x)-q_i)y_i=\sum_i x_i (G_i(y)-q^\prime_i),$$ for all $x,y \in X$, 
  where $q_i=F_i(0), q^\prime_i=G_i(0)$. Then $F$ (hence also $G$) is affine, i.e. there are $a_{ij} \in \clc^{\rm s.a.}$ such that $F_i(x)=q_i+\sum_j x_j a_{ji}$ for all $i$.
\elmma
 {\it Proof:}\\
 The idea of the proof is quite similar to the Lemma \ref{iso_aff} above. First by replacing $F_i, G_i$ by $F_i(\cdot)-q_i$ and $G_i(\cdot)-q^\prime_i$ respectively, 
  we may assume without loss of generality that $q_i=q^\prime_i=0$. Then  observe that, if $\sum_l c_l x^{(l)} =0$ for $x^{(1)}, \ldots, x^{(p)} \in X$ and 
  $c_1, \ldots, c_p \in \IR$, we have for any $y \in X$ 
   $\sum_l c_l (\sum_{i=1}^n F_i(x^{(l)})y_i=\sum_{i=1}^n \sum_l c_l x^{(l)}_i G_i(y)=0$. That is, $\sum_i Q_i y_i=0$ for all $y \in X$, where $Q_i=\sum_l c_l F_i(x^{(l)}).$ But as 
    the coordinate functions restricted to $X$ form a linearly independent set, we conclude $Q_i=0$ for each $i$. This implies that there is a well-defined linear extension of 
     $F$ on $\IR^n$, as in Lemma \ref{iso_aff}, and the elements $a_{ij}$ can be obtained in a similar way. \qed\\

   \subsection{Main result : existence of a universal quantum group of isometry}
       Throughout this section, let $X \subset \IR^n$ be a compact subset with $d$ the restriction of Euclidean metric inherited from $\IR^n$, i.e. $d^2(x,y)=\sum_{i=1}^n (x_i-y_i)^2$.
     Let $X_1, \ldots, X_n$ denote the restriction of the coordinate functions of $\IR^n$ to the subset $X$. Thus, $C(X)$ is generated as a $C^*$ algebra by $1$ and $X_1, \ldots, X_n$. 
 Without loss of generality (if necessary by translating the set $X$) we can assume that $0$ of $\IR^n$ belongs to $X$ and moreover, 
$X_1, \ldots, X_n$ are linearly independent. This also implies the linear independence of $\{ 1, X_1, \ldots, X_n \}$. Indeed, if $X_i$'s are linearly dependent,
say $X_1, \ldots, X_k$ are independent and $X_j=\sum_{l=1}^k d_{jl}X_l$ for $j=k+1, \ldots, n$, we can write the metric $d^2(x,y)$ as 
$\underline{Z}^\prime (I_k+D^\prime D) \underline{Z}$, where $\underline{Z}=((x_1-y_1),\ldots, (x_k-y_k))$ and $D=(( d_{jl}))$. 
That is, considering new coordinate functions $\hat{X}_i= \sum_{j=1}^k c_{ij}X_j$, with $C=(( c_{ij}))$ where $C=(I+D^\prime D)^{\frac{1}{2}}$ which is 
positive and invertible, we have linearly independent $\hat{X}_i$ ($i=1,\ldots,k$) such that $d^2=\sum_{i=1}^k (\hat{X}_i \ot 1-1 \ot \hat{X}_i)^2.$ 

Let $\alpha$ be a faithful action of a CQG $\clq$ on $C(X)$ and let $\alpha_r$ be the reduced action of the reduced 
 CQG $\clq_r$ with bounded antipode say $\kappa$, $\clc_0$ be the spectral subalgebra of $C(X)$ corresponding to $\alpha$, $F_i(\cdot):=\alpha(X_i)(\cdot),$ $F^r_i(\cdot)=\alpha_r(X_i)(\cdot).$ 
  We first want to derive several  equivalent conditions for the action $\alpha$ to be isometric. These will be very similar to what one gets for classical groups acting isometrically 
    on subsets of Euclidean spaces.

\bthm
\label{iso_equiv}
The following are equivalent.\\
(i) $\alpha$ is isometric.\\
(ii) \be \label{666} \sum_i (x_i1-F^r_i(y))^2=\sum_i (\kappa(F^r_i(x)-y_i 1))^2\ee $ \forall x,y \in X$.\\
(iii) $F_i \in \clc_0 \ot_{\rm alg} \clq_0$  $\forall i$ and  $$ \sum_i (x_i1-F_i(y))^2=\sum_i (\kappa(F_i(x)-y_i 1))^2$$ $ \forall x,y \in X$.\\
(iv)  $F_i \in \clc_0 \ot_{\rm alg} \clq_0$  $\forall i$ and 
\be \label{rst} \sum_i (F^2_i(x)+F^2_i(y)-2F_i(x)F_i(y))=d^2(x,y)1,\ee
$\forall x,y \in X.$\\
(v)  $F_i \in \clc_0 \ot_{\rm alg} \clq_0$  $\forall i$ and 
\be \label{pqr}  \sum_{i=1}^n (F_i(x)-F_i(y))^2=d^2(x,y)1.\ee
\ethm
{\it Proof:}\\
We begin with the observation that (ii) is nothing but the equivalent condition of isometry (\ref{iso_cond}) obtained in 
 Lemma \ref{iso_cond_lemma}. Clearly, (iii) implies (ii). To see the implication $(ii) \Longrightarrow (iii)$ write $H_i=F^r_i$, $q_i=H_i(0)$, $ G_i(x)=\kappa(H_i(x)),$
  $q^\prime_i=\kappa(q_i)$. Putting $x=0$ in (\ref{666}) we get \be \label{777} \sum_i H_i(y)^2=\sum_i (y_i1-q^\prime_i)^2,\ee which also gives, by applying $\kappa$ and using $\kappa^2={\rm id}$,
  \be \label{898} \sum_i G_i(y)^2=\sum_i (y_i1-q_i)^2,\ee  and $\sum_i q_i^2=\sum_i (q^\prime_i)^2$ by 
   putting $y=0$. Expanding (\ref{666}) and using (\ref{777}), (\ref{898}) as well as $\sum_i q_i^2=\sum_i (q^\prime_i)^2$ we get, on simplification, 
   $$ \sum_i x_i (H_i(y)-q_i)=\sum_i y_i(G_i(x)-q^\prime_i).$$ By Lemma \ref{affine_action} we conclude that $H_i$, i.e. $F^r_i$ is affine, 
    say, $F^r_i=1 \ot q_i +\sum_j X_j \ot a_{ji}$ for $q_i, a_{ji} \in \clq^{s.a.}_r$. But as noted before (Proposition 3.8, \cite{alex_preprint}), 
     $\alpha_r$ is injective and induces a unitary representation. Clearly, $q_i, a_{ij}$'s belong to the span of 
      the matrix coefficients of the restriction of the above unitary representation to the $n+1$ dimensional invariant 
       subspace spanned by $\{ 1,X_1, \ldots, X_n\}$, hence $q_i,a_{ji} \in \clq_0$, i.e.  $F^r_i=\alpha_r(X_i) \in C(X) \ot_{\rm alg} \clq_0$. 
        Applying Proposition 2.2 of \cite{soltan} to the injective action $\alpha_r$, we conclude  that $X_i \in \clc_0$ and thus $F^r_i \in \clc_0 \ot \clq_0$, 
        which also means $F_i=F^r_i \in \clc_0 \ot_{\rm alg} \clq_0$ as $\pi_r={\rm id}$ on $\clq_0$.\\

$(iii) \Longrightarrow (iv):$\\
Use Sweedler notation $F_i=X_{i(0)} \ot X_{i(1)}$. Using $\kappa^2={\rm id}$ we get from (iii) the following:\\
\bea \label{888} \lefteqn{\sum_i \left( X_i^2 \ot 1 \ot 1+1\ot X_{i(0)}^2 \ot \kappa(X_{i(1)}^2)-2 X_i \ot X_{i(0)} \ot \kappa (X_{i(1)}) \right)}{\nonumber}\\
&=& \sum_i 
\left( X_{i(0)}^2 \ot 1 \ot X_{i(1)}^2 +1 \ot X_i^2 \ot 1-2 X_{i(0)} \ot X_i \ot X_{i(1)} \right).\eea
 Consider  $\Psi:=({\rm id} \ot m) \circ (\alpha \ot {\rm id})$  from $ \clc_0 \ot_{\rm alg} \clq_0$ to itself 
  where $m$ denotes the multiplication of $\clq_0$.  Using $({\rm id} \ot \epsilon)(\alpha(f))=f$ as well as
  $m \circ ({\rm id} \ot \kappa) \circ \Delta(a)=\epsilon(a) 1$ for all $a \in \clq_0$, we get for $f \in \clc_0$ the following:
    $\Psi(f_{(0)} \ot \kappa(f_{(1)}))=f_{(00)} \ot f_{(01)}\kappa(f_{(1)})=f_{(0)} \ot f_{(11)}\kappa(f_{(12)})=f_{(0)}\epsilon(f_{(1)})\ot 1=f \ot 1$.
     Thus, by applying $( {\rm id} \ot \Psi)$ on both sides of (\ref{888}) we get
   $$ \sum_i (F_i(x)^2+F_i^2(y)-2F_i(y)F_i(x)=d^2(x,y)1.$$ 

 $(iv) \Longrightarrow (v):$\\
 Interchanging 
  $x,y$, using $d(x,y)=d(y,x)$ and adding we get $$ \sum_i \left(F_i^2(x)+F_i^2(y)-F_i(x)F_i(y)-F_i(y)F_i(x)\right)=d^2(x,y)1, $$ or, equivalently, 
  $$ \sum_i (F_i(x)-F_i(y))^2= d^2(x,y)1 \equiv \sum_i (x_i-y_i)^2 1.$$
  
  $(v) \Longrightarrow (iii):$\\
 It follows by applying $\sigma_{23} \circ (\Psi \ot {\rm id})$ (where $\sigma_{23}$ flips the second and third tensor
  components) on both sides of (\ref{pqr}).
    \qed\\
  
  \bcrlre
  \label{cor_aff}
  Assume the set-up of Theorem \ref{iso_equiv}. Then $\alpha$ is isometric if and only if the following hold:\\
   (a) $F_i\equiv \alpha(X_i)=\sum_j X_j \ot a_{ji} + 1\ot \xi_i$ for some  self-adjoint elements $a_{ij}, \xi \in \clq_0$.\\
(b) \be \label{1000} \sum_i a_{ji}a_{ki} = \sum_i a_{ij} a_{ik} =\delta_{jk} 1, \ee for all $j,k$ where $\delta_{jk}$ is the Kronecker's delta.\\
  (c) There is an $\alpha$-invariant probability measure, say $\mu$, on $X$ such that 
   $\xi_i=c_i 1-\sum_j c_j a_{ji}$ for each $i$, where $c_i=\int_X X_i d \mu$, hence 
  $\alpha$ satisfies 
  \be \label{8888} \alpha(X_i-c_i1)=\sum_j (X_j-c_j1) \ot a_{ji}, \ee for all $i=1,\ldots,n.$\\
  (d) There is a surjective CQG morphism $\phi : A_o(n) \raro \clq$ which sends the canonical generators $q^o_{ij}$ of $A_o(n)$ to $a_{ij}$.
    \ecrlre
  {\it Proof:}\\
It is easy to check from the defining relations of the generators of $A_o(n)$  that (a)-(d) imply (\ref{pqr}), hence isometry of the action. 

 Conversely, assuming $\alpha$ to be isometric,  we already obtained (a) in the proof of Theorem \ref{iso_equiv}. 
  To prove (b) we first use $\sum_i (F_i(x)-F_i(y))^2=\sum_i (x_i-y_i)^21$ for all $x,y \in X$, with $y=0$ to get 
 \be \label{2000} \sum_i \left( \sum_j x_ja_{ji} \right)^2 =\sum_i x_i^2 1.\ee From the isometry condition we get 
 $$ \sum_i F_i(x)^2 +\sum_i F_i(y)^2 -2 \sum_i F_i(x)F_i(y)=(\sum_i x_i^2+\sum_i y_i^2 -2 x_i y_i )1.$$ 
 Using (\ref{2000}) in this and by simplification, we get for all $x,y \in X$ \be \label{3000} \sum_j (x_j-y_j) \left( \sum_i (\xi_i a_{ji}-a_{ji} \xi_i) \right) -2 \sum_{j,k} x_jy_k \left( \sum_i a_{ji}a_{ki} \right)=-2 (\sum_j x_j y_j ) 1.\ee Putting $y=0$ in this and using linear independence of the coordinate functions $X_j$'s we conclude $\sum_i (\xi_i a_{ji}-a_{ji} \xi_i)=0$ for all $j$, and plugging it  back into (\ref{3000}) we have
  $$ \sum_{j,k} x_jy_k (\sum_i a_{ji}a_{ki})=(\sum_j x_jy_j)1.$$ Using linear independence of the $x_i$'s, with fixed value of $y$ and each $j$, 
  we get  $ \sum_k y_k (\sum_i a_{ji} a_{ki})=y_j 1.$ As this is true for all $y \in X$ again 
  appealing to the linear independence of the coordinate functions we conclude \be \label{4000} \sum_i a_{ji}a_{ki}=\delta_{jk}1.\ee
 
 Now, choose any faithful state $\phi$ on $C(X)$ and consider the $\alpha$-invariant state $\overline{\phi}=(\phi \ot h)\circ \alpha$ where $h$ is the Haar state. 
  Let $\mu$ be the corresponding probability measure on $X$. Clearly, faithfulness of $h$  on $\clq_0$ implies faithfulness of $\overline{\phi}$  on 
   $\clc_0$, hence in particular on the unital $\ast$-subalgebra generated by $X_i$'s, as we observed in the proof of 
    Theorem \ref{iso_equiv} that $X_i \in \clc_0$ $\forall i$.
   
    Now, letting $c_i=\int_X X_i d \mu=\overline{\phi}(X_i)$, we have from the $\alpha$-invariance 
  the following: $$ c_i 1=\xi+\sum_j c_j a_{ji},$$ i.e. $\xi_i=c_i1-\sum_j c_j a_{ji}.$ From this, we get 
  \be \label{333} \alpha(X_i-c_i1))=\sum_j (X_j-c_j1) \ot a_{ji}.\ee 
  This proves (c). Moreover, the span of $Y_1=X_1-c_11, \ldots, Y_n=X_n-c_n1$, say $\clv$, is a finite-dimensional subspace of the spectral subalgebra 
   $\clc_0$, and (\ref{333}) implies that it is left invariant by $\alpha$. Therefore,  $\alpha|_\clv$ must be a nondegenerate finite-dimensional representation of $\clq$ and 
     the corresponding $\clq$-valued matrix w.r.t. the basis $Y_1,\ldots, Y_n$, i.e. $(( a_{ji} ))$, must be invertible in $M_n(\clq)$. But (\ref{4000})
  already shows the one-sided inverse to be $((a_{ij} ))$, so this must be the both-sided inverse as well, giving $\sum_j a_{ij} a_{ik}=\delta_{jk} 1$ and 
   thereby also completing the proof of (b).
   
  The proof the statement (d) follows from the universality of $A_o(n)$.
  \qed\\   
 
 The next step is to show that the constants $(c_1, \ldots, c_n)$ in (c) of Corollary \ref{cor_aff} can be chosen independently of the CQG or the action. In other words, there are 
  universal such constants depending only $(X,d)$ which work for every isometric CQG action. 
\blmma
\label{univ_const}
Let $X \subseteq \IR^n$ be a compact subset and let $d$ be the restriction of the Euclidean metric of $\IR^n$ to $X$. Then there are real numbers $(c_1, \ldots, c_n)$ depending only on 
 $(X,d)$ such that for every CQG $\clg$ with a faithful isometric action $\alpha_\clg$ of $\clg$ on $(X,d)$, there is a surjective CQG morphism $\pi_\clg: A_o(n) \raro \clg$ such that 
  $$ \alpha_\clg(X_i-c_i1)=\sum_j (X_j -c_j1) \ot g_{ji},$$ where $g_{ji}=\pi_\clg(q^o_{ji})$, $q^o_{ji}$ being the canonical generators of $A_o(n)$ discussed earlier.
\elmma
{\it Proof:}\\
Let $X_1,\ldots, X_n$ be the coordinate functions of $X$ as in the Lemma \ref{iso_equiv}. 
Let $M$ be a large enough positive number such that $X \subseteq [-M, M]^n$, so that each $X_i$ and hence also $\mu(X_i)$ for 
 any probability measure on $X$ is bounded above by $M$ and let $M^\prime=(n+1)M$. Let $\clb$ the free product of $n$ copies of $C([-M^\prime,M^\prime])$ and $\cls=A_o(n) \ast \clb$. Denote the canonical generators 
  of $A_o(n)$ by $q^o_{ij}$ and the coordinate function of the $i$-th copy of $C[-M^\prime,M^\prime]$ 
  in the free product $\clb$ by $\beta_i$. Thus, $\{ q^o_{ij}, \beta_i, i,j=1,\ldots, n\}$ generate 
   $\cls$ as a $C^*$ algebra.

         We have noted in Corollary \ref{cor_aff}  that given any faithful isometric action $\alpha_\clg$ of a CQG $\clg$ on $(X,d)$,
          there are  generating elements $g_{ij}$ for the CQG $\clg$ and a surjective morphism
          $\pi_\clg$ from $A_o(n)$ to $\clg$ which sends $q^o_{ij}$ to $g_{ij}$. Moreover,  there are  $c^g_i=\int_X X_i d \mu$ for some invariant probability 
           measure on $X$ such that $\alpha_\clg (X_i-c^g_i1)=\sum_j (X_j -c^g_j1) \ot g_{ji}.$ It is clear that $|c^g_i| \leq M$ and hence $\| c^g_i1-\sum_j c^g_j g_{ji} \|
            \leq M^{\prime }$. Define, by the universality of free product, a surjective $C^*$ homomorphism from $\cls$ to $\clg$ sending $q^o_{ij}$ to $g_{ij}$ and $\beta_i$ to $c^g_i1-
             \sum_j c^g_j g_{ji}$. Call this $\rho_\clg$ and let $\cli_\clg$ be the kernel of it, which is a closed two-sided ideal.
         Let us now consider the collection $\clf$  of    all  $\cli_\clg$ corresponding to  faithful isometric CQG actions  $(\clg, \alpha_\clg)$    
           and let $\cli^0$ be the intersection of all such $\cli_\clg$'s. Let $\clq=\cls / {\cli^0}$, $b_{ij}=q^o_{ij}+\cli^0,$ $\gamma_i=\beta_i+\cli^0$ 
            in $\clq$. Denote by $\theta_\clg$ the quotient map from $\clq=\cls/\cli$ to $\cls/ {\cli_\clg}$. 
            To prove the lemma it suffices to show that there is some $c \equiv (c_1,\ldots, c_n) \in [-M^\prime, M^\prime]^n$ as above such that $\gamma_i=c_i1_\clq-\sum_j c_j b_{ji}$.

            To this end,  consider the Banach space $\clq^{(n)}$ which is the direct sum of $n$ copies of $\clq$. Let $\xi^{(i)},i=1, \ldots, n$ 
             be elements of $\clq^{(n)}$ defined by $\xi^{(i)}=(\xi^{(i)}_1,\ldots, \xi^{(i)}_n)$, where $\xi^{(j)}_i=-b_{ji}$ if $j \ne i$ and 
             $\xi^{(i)}_i=1-b_{ii}$. Let $\gamma:=(\gamma_1, \ldots, \gamma_n) \in \clq^{(n)}$. 
             
           Consider the finite dimensional subspace  $\clk$ of the Banach space $\clq^{(n)}$ spanned by   $\xi^{(i)}, i=1,\ldots, n$ and let $\clv$ 
            be the  (finite dimensional) subspace spanned by $\clk$ and $\gamma$. We want to prove that $\gamma \in \clk$.
              Consider the subspace $\clw$ of the dual of $\clq^{(n)}$ spanned by 
             the bounded linear functionals of the form $ (\psi_1 \circ \pi_{\clg},\ldots \psi_n \circ \pi_\clg) $  where $\cli_\clg \in \clf$, $\psi_i$ is any  bounded linear functional on $\cls/{\cli_\clg}$
              and $\pi_\clg : \cls \raro \cls/{ \cli_\clg}$ is 
              the quotient map. It is easy to see that for some $q=(q_1, \ldots, q_n)\in \clq^{(n)}$, $w(q)=0$ for all $w \in \clw$ implies 
               $q_i+\cli_\clg =0$ for all $i$ and all $\cli_\clg$, hence $q_i \in \cli^0$ for all $i$, i.e. $q=0$ as an element of $\clq^{(n)}$. 
               This means $\clw$ is weak-$\ast$ dense in the dual of $\clq^{(n)}$, so in particular $\{ w|_\clv:~w \in \clw \}$ is dense in the dual of the finite 
               dimensional space $\clv$, hence must coincide with it. Thus, if $\gamma$ is not in $\clk$, we can get $w \in \clw$, say of the form $(w_1 \circ \pi_\clg, \ldots, 
               w_n \circ \pi_\clg)$ for some $\clg$ and $w_i \in \clq^*$,   such that  $w(\gamma)$ is nonzero and $w$ vanishes on $\clk$.  But this is
                not possible by definition of $\cli_\clg$,  because we have $(c_1^g, \ldots c^g_n)$ such that $(\gamma_i+\cli_\clg)=\sum_j c^g_j (\xi^{(j)}+\cli_\clg)$ for all $i$, 
                which implies $w(\gamma)=\sum_j c^g_j w(\xi^{(j)})=0$ This completes the proof. \qed\\

    Now we are in a position to prove the existence of the universal object in the category of all CQG acting isometrically.

  \bthm
\label{existence_univ}
Let $X$ be a compact subset  $\IR^n$ and $d$ be the restriction of the Euclidean metric. 
Then  the category ${\bf Q}(X,d)$   has a universal object to be denoted by $QISO(X,d)$.
\ethm 
{\it Proof:}\\ We use the notation in the statement as well as proof of the previous Lemma \ref{univ_const}. 
 We claim that $\clq$ is actually a CQG and it is indeed the desired universal object. For this,  
                       let us now consider a category ${\bf  C}$  with objects $(\clc, \alpha_\clc, \{ x_{ij},i,j=1, \ldots, n \})$ where 
                        $\clc$ is a unital $C^*$ algebra  generated by self-adjoint elements 
                        $\{ x_{ij} , i,j=1, \ldots, n\}$ such that $(( x_{ij} ))$ as well as $(( x_{ji} ))$ are  unitaries and also there is a unital $\ast$-homomorphism $\alpha_\clc$ 
                         from $C(X)$ to $C(X) \ot \clc$ sending $(X_i-c_i1)$ to $\sum_j (X_j -c_j1) \ot x_{ji}$. The morphisms from $(\clc, \alpha_\clc, \{ x_{ij}\})$
                          to $(\cld, \alpha_\cld, \{ y_{ij} \} )$  are unital $\ast$-homomorphisms $\beta : \clc \raro \cld$ such that 
                          $\pi(x_{ij})=y_{ij}$ for all $i,j$. This object-class of this category  is clearly nonempty and  it  has $(\clq, \alpha_{\cli^0}, \{ b_{ij} \} )$
                           in it. 
                          
                          Moreover, by definition 
                          of each object $(\clc, \alpha_\clc, \{ x_{ij} \})$  we get a unital surjective $\ast$-homomorphism (say $\rho_\clc$) from $A_o(n)$ to 
                          $\clc$ sending $q^o_{ij}$ to $x_{ij}$. Let the kernel of this 
                           map be $\cli_\clc$ and let $\cli$ be the intersection of all such ideals and let $\clm:=A_o(n)/\cli$. We claim that $\clm$ is the universal object 
                            in ${\bf C}$ and it is also a CQG. Denote the collection of all $\clc$ corresponding to objects of ${\bf C}$ by $\clf$ and the 
                             quotient map from $A_o(n)$ to $\clm$ by $\pi$. Moreover, let $\eta_\clc$ be the canonical map from $\clm$ onto $\clc$ such that $\eta_\clc \circ \pi=\rho_\clc$,
                              and we have $\bigcap_{\clc \in \clf} {\rm ker}(\eta_\clc)=(0),$ as $\cli={\rm Ker}(\pi)=\bigcap_\clf \cli_\clc$. 
                              
                              Clearly, $\clm$ is generated by $m_{ij}=\pi(q^o_{ij}), i,j=1, \ldots, n$.
                            It suffices to show that there is an action $\alpha: C(X) \raro C(X) \ot \clm$ satisfying $\alpha(X_i-c_i1)=\sum_j (X_j-c_j1) \ot m_{ji}$.
                             To this end, we choose sufficiently large $M>0$ such that $|x_i-c_i|\leq M$ for all $i$ and all $x=(x_1,\ldots, x_n) \in X$. Let $M^\prime=nM$. Consider  
                               $F_i(x):=c_i1+\sum_j(x_j-c_j)q^o_{ji} \in A_o(n)$, $x=(x_1,\ldots, x_n) \in X$, $i=1, \ldots, n$. Now,  $\forall i,j$, $\forall x \in X$ and 
                                $\clc \in \clf$, the fact that $\alpha_\clc$ is a homomorphism implies
                                 $\rho_\clc([F_i(x), F_j(x)]) =0$, i.e. $[F_i(x), F_j(x)] \in \cli_\clc$ $\forall \clc \in \clf$, hence $[F_i(x), F_j(x)] \in \cli$, i.e.
                                 $\pi([F_i(x), F_j(x)])=[ \pi(F_i(x)), \pi(F_j(x)) ]=0$. This means $F^\prime_i \in C(X) \ot \clm$ ($i=1, \ldots, n$)
                                  given by $F^\prime_i(x):=\pi(F_i(x))$ are self-adjoint 
                                  mutually commuting elements.  Moreover, as $\| m_{ij}\|=\|  \pi(q^o_{ij}) \| \leq 1$ for all $i,j$ and $|x_i-c_i| \leq M$
                                  for all $x=(x_1, \ldots, x_n) \in X$,
                                   we have $\| F_i^\prime-c_i 1 \ot 1 \| \leq nM$. Consider the set $Y(\subseteq \IR^n ):=\{ (y_1, \ldots, y_n):~ |y_i-c_i| \leq nM~\forall i\}.$ 
                                    Denote by $X^\prime_1, \ldots, X^\prime_n$ the restrictions of coordinate functions to $Y$.   Then $C(Y)$ is 
                                    the  unital $C^*$ algebra generated by $n$ mutually commuting self-adjoint elements $X^\prime_1, \ldots, X^\prime_n$ satisfying 
                                     $\| X^\prime_i-c_i1\| \leq nM$, and it is in fact the universal 
                                     such $C^*$ algebra. Thus, we get a $C^*$ homomorphism, say $\tilde{\alpha}$, from $C(Y)$ to $C(X) \ot \clm$ which sends $X^\prime_i$ to $F^\prime_i$ for 
                                      all $i$. Clearly, by construction, $({\rm id} \ot \eta_\clc) \circ \tilde{\alpha}(X^\prime_i)=\alpha_\clc(X_i)$ (where $\eta_\clc: \clm \raro \clc$
                                       as before and $X_i$ is the restriction 
                                       of the $i$-th coordinate function to $X \subset Y$, i.e. $X_i=X^\prime_i|_X$) for all $\clc \in \clf$, hence 
                                       by continuity and density of polynomials in $C(Y)$, we have $({\rm id} \ot \eta_\clc) \circ \tilde{\alpha}(f)=\alpha_\clc(f|_X)$ for all $f \in C(Y)$.
                                        In particular, if $f|_X=0$, we get $\eta_\clc(\tilde{\alpha}(f)(x)) =0$ for all $x \in X$ and $\clc \in \clf$, hence $\tilde{\alpha}(f) =0$ in
                                         $C(X) \ot \clm$ as $\bigcap_{\clc \in \clf} {\rm Ker}(\eta_\clc)=(0)$. In other words, the ideal of $C(Y)$ consisting of all $f$ which is identically $0$ on the subset $X$ is in the kernel of $\tilde{\alpha}$, 
                                          so $\tilde{\alpha}$ descends to the quotient of $C(Y)$ by this ideal, which is nothing but $C(X)$. Denoting this induced map by $\alpha$, we do get a 
                                           $C^*$ homomorphism from $C(X)$ to $C(X) \ot \clm$ satisfying $\alpha(X_i-c_i1)=\sum_j (X_j-c_j1) \ot m_{ji}$. This gives us an object $(\clm, \alpha, \{ m_{ij} \} )$
                                            in ${\bf C}$, which is clearly universal by construction.

                           To obtain the coproduct on $\clm$, consider  $M_{ij}:=\sum_k m_{ik} \ot m_{kj} \in \clm \ot \clm$ and note 
                           that $M_{ij}=(\pi \ot \pi )(Q_{ij}),$ where $Q_{ij}=\sum_k q^o_{ik} \ot q^o_{kj},$ $i,j=1,\ldots, n$. From the definition of $A_o(n)$ one can check that 
                            the unital $C^*$-subalgebra (say $\cld$) of $\clm \ot \clm$ generated by $Q_{ij}$'s along with the map 
                            $\beta=(\alpha \ot {\rm id}) \circ \alpha : C(X) \raro C(X) \ot \cld$ give an object in ${\bf C}$, hence 
                             by universality of $\clm$ there is a well-defined unital $\ast$-homomorphism $\Delta_0$ from 
                             $\clm$ to $\cld \subseteq \clm \ot \clm$. Clearly, $\Delta_0$ satisfies $\Delta_0 \circ \pi=(\pi \ot \pi) \Delta$, where 
                             $\Delta$ is the coproduct of the CQG $A_o(n)$. It follows that $\Delta$ maps ${\rm ker}(\pi)$ to ${\rm Ker}(\pi \ot \pi)$ and 
                             moreover, from the density of each of the linear spans of $\Delta(\clq_0) (1 \ot \clq_0)$ as well as $\Delta(\clq_0)(\clq_0 \ot 1)$ in $A_o(n) \ot A_o(n)$ where $\clq_0$ is the $\ast$-algebra generated by 
                              $q^o_{ij}$'s, we get (by applying $\pi \ot \pi)$ similar density with $A_o(n)$ and $\Delta$ replaced by $\clm$ and $\Delta_0$ respectively,
                               and $\clq_0$ by the algebra generated by the $m_{ij}$'s. Therefore, $\cli={\rm Ker}(\pi)$ is a closed Hopf ideal and 
                              $\clm$ becomes a quantum subgroup of $A_o(n)$, hence a CQG in particular, with the coproduct $\Delta_0$. 
                             The map $\alpha$ now becomes an action of the CQG $\clm$.
                              Clearly, $\alpha$ is 
                               an isometric faithful action (by Corollary \ref{cor_aff}),
                                hence $\clm$ must be a quotient of $\clq$. But on the other hand, by the universality of $\clm$ and the fact that $\clq$ is an object of 
                               ${\bf C}$, it would follow that $\clq\cong \clm$. 
                                                                                              \qed\\

 
The following corollary shows that we can relax the assumption of isometric embedding of $(X,d)$ into Euclidean space to some extent, 
 and it actually suffices  to have a topological embedding of $X$ into some $\IR^n$ so that the metric inherited from the Euclidean metric on the image of $X$ is bijectively 
  related to the original metric $d$. 

  \bcrlre
  \label{3334}
  Let $(X,d)$ be a compact metric space. Suppose also that  there are topological embedding $f: X \raro \IR^n$ and
    a homeomorphism $\psi$ of $\IR^+$  such that  
   $(\psi \circ d)(x, y)=d_0(f(x),f(y))$ for
    all $x,y \in X$, where we have denoted the Euclidean metric of $\IR^n$ by $d_0$. 
    Then the conclusion of Theorem \ref{existence_univ} holds.
    \ecrlre
    {\it Proof:}\\
    It suffices to note that a CQG action on $X$ is isometric w.r.t.  $d$ if and only if it is isometric w.r.t. the metric $d_0$ on $X \cong f(X) \subset \IR^n$. \qed\\

    \brmrk
    It follows from \cite{embedding_thm} that an arbitrary finite metric space satisfies 
     the condition of Corollary \ref{3334} with $\psi(t)=t^c$ for some $c>0$. This implies that our existence theorem does extend that of Banica for finite spaces. 
     Examples of  metric spaces satisfying the condition of the Corollary \ref{3334} also include the spheres $S^n$ for all $n \geq 1$.
 \ermrk
 If the metric space $X$ in Corollary \ref{3334} has at least $4$  components each of which is isometric to some given set,   $QISO(X,d)$ will have $\cls_4^+$ as a quantum subgroup, 
 hence genuine.
  It is more interesting   to 
   construct examples of genuine CQG acting isometrically on   connected spaces. We refer to  \cite{huang1} (see also \cite{etingof}) for a rich source of such examples and Section 5 
    of the present article of computations of  quantum isometry groups of some of these spaces. 
   In fact, for a compact connected $X \subseteq \IR^n$, $QISO(X,d)$ coincides with $C(ISO(X,d))$ in the following two cases:\\
   (i) $X$ is an embedded submanifold,\\
   (ii) $X$ has nonempty interior in $\IR^n$.\\
  Both these statements follow from \cite{nonexistence} because any isometric action is affine. This implies  smoothness in case (i) and we can apply Theorem 10.6 of 
  \cite{nonexistence}.  
    In case (ii), the conclusion follows from  Lemma 10.1 of \cite{nonexistence}.

For a compact connected Riemannian manifold  $M$ one has  a natural metric $d$ coming from the geodesic distance and it is interesting to compare 
 the corresponding quantum isometry group $QISO(M,d)$ (whenever it exists) with the geometric quantum isometry group (say $QISO^\cll(M)$) defined in \cite{goswami} in terms of the 
  Hodge-de-Rham Laplacian. In fact, it follows from  \cite{nonexistence} that $QISO^\cll(M)=C(ISO(M))$ in this case and as we have already noted, whenever $M$ satisfies the hypothesis of 
   Corollary \ref{3334},  $QISO(M,d)$ exists and equals $QISO^\cll(M)=C(ISO(M)).$  We believe that it is true in general. A. L. Chirvasitu seems to have a proof 
    of this when $M$ is negatively curved (private communication).

 \section{Computation of $QISO(X,d) $ for a class of metric spaces}
 In this final section
   we compute quantum isometry groups of  two  classes of compact spaces $T_n$ and $T^\prime_n$, $n=1,2, \ldots$ shown in the figures below.
 These are particular types of examples considered in \cite{huang1} (see also \cite{etingof}).

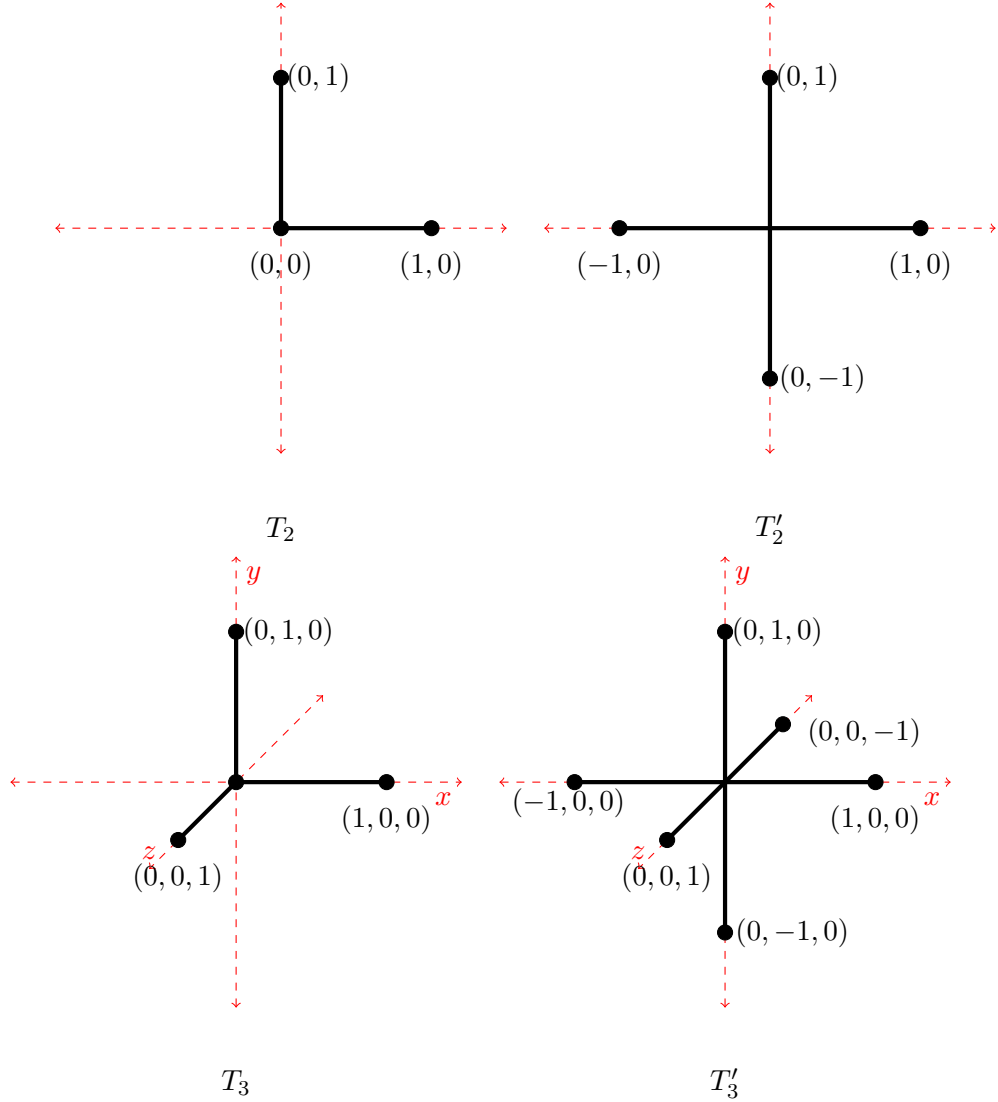
\begin{tikzpicture}
  
\draw[dashed,red, <->] (-3,0) -- (3,0);
\draw[dashed, red, <->] (0,-3) -- (0,3) ;
\draw[ultra thick, -](0,0)--(2,0);
\draw[ultra thick,-](0,0)--(0,2);
\draw [fill] (0,0) circle [radius=0.1];
\draw[fill](0,2) circle [radius=0.1];
\draw[fill](2,0) circle [radius=0.1];
\node at (0,-4){ $ T_2$};
\node at (-0,-0.5) {$(0,0)$};
\node at (2,-0.5){$(1,0)$};
\node at (0.5, 2){$(0,1)$};
\begin{scope}[xshift=6.5cm]
\draw[dashed,red, <->] (-3,0) -- (3,0);
\draw[dashed, red, <->] (0,-3) -- (0,3) ;
\draw[ultra thick, -](-2,0)--(2,0);
\draw[ultra thick,-](0,-2)--(0,2);
\draw [fill] (-2,0) circle [radius=0.1];
\draw[fill](2,0) circle [radius=0.1];
\draw [fill] (0,-2) circle [radius=0.1];
\draw[fill](0,2) circle [radius=0.1];
\node at (-2,-0.5) {$(-1,0)$};
\node at (2,-0.5){$(1,0)$};
\node at (0.7, -2){$(0,-1)$};
\node at (0.5, 2){$(0,1)$};
\node at (0,-4){ $ T_2^\prime$};
\end{scope}
 \end{tikzpicture}
  \begin{tikzpicture}
    \draw[dashed,red, <->] (-3,0,0) -- (3,0,0) node[anchor=north east]{$x$};
\draw[dashed,red,<->] (0,-3,0) -- (0,3,0) node[anchor=north west]{$y$};
\draw[dashed,red, <->] (0,0,-3) -- (0,0,3) node[anchor=south]{$z$};
\draw[ultra thick, -](0,0,0)--(2,0,0);
\draw[ultra thick,-](0,0,0)--(0,2,0);
\draw[ultra thick, -](0,0,0)--(0,0,2);
\draw [fill] (0,0,0) circle [radius=0.1];
\draw[fill](0,2,0) circle [radius=0.1];
\draw[fill](2,0,0) circle [radius=0.1];
\draw[fill](0,0,2) circle [radius=0.1];
\node at (2,-0.5,0){$(1,0,0)$};
\node at (0.7, 2,0){$(0,1,0)$};
\node at (0,-0.5,2){$(0,0,1)$};
\node at (0,-4,0){ $ T_3$};
 \begin{scope}[xshift=6.5cm]
  \draw[dashed,red, <->] (-3,0,0) -- (3,0,0) node[anchor=north east]{$x$};
\draw[dashed,red,<->] (0,-3,0) -- (0,3,0) node[anchor=north west]{$y$};
\draw[dashed,red, <->] (0,0,-3) -- (0,0,3) node[anchor=south]{$z$};
\draw[ultra thick, -](-2,0,0)--(2,0,0);
\draw[ultra thick,-](0,-2,0)--(0,2,0);
\draw[ultra thick, -](0,0,-2)--(0,0,2);
\draw[fill](0,2,0) circle [radius=0.1];
\draw[fill](2,0,0) circle [radius=0.1];
\draw[fill](0,0,2) circle [radius=0.1];
\draw[fill](0,-2,0) circle [radius=0.1];
\draw[fill](-2,0,0) circle [radius=0.1];
\draw[fill](0,0,-2) circle [radius=0.1];
\node at (2,-0.5,0){$(1,0,0)$};
\node at (0.7, 2,0){$(0,1,0)$};
\node at (0,-0.5,2){$(0,0,1)$};
\node at (-2,-0.2,0.2){$(-1,0,0)$};
\node at (0.9, -2,0){$(0,-1,0)$};
\node at (0.7,-0.5,-3){$(0,0,-1)$};
\node at (0,-4,0){ $ T_3^\prime$};
 \end{scope}
 \end{tikzpicture}

 Formally,  $ T_n$ (respectively $T^\prime_n$) is a subset of $\IR^n$ obtained by gluing $n$ copies of an interval $I$ (which is $[0,1]$ or $[-1,1]$ respectively) 
 at the origin ${\bf 0}$. The $i$-th copy of $I$, denoted by  $T_i$ or $T^\prime_i$ respectively,  is identified with the subset $(0,0,\ldots, 0,t,0,\ldots, 0),~t \in I$ of 
  $\IR^n$, where $t$ is at the $i$-th place.
  
  Here, $d$ will denote the restriction of the Euclidean metric. Let  $X_i$ be the restriction of the $i$-th coordinate function to $T_n$ or $T^\prime_n$ and let 
   $\clq$ be the corresponding quantum isometry group whose (affine, faithful) action is given by $$\alpha(X_i)=\sum_{j=1}^n X_j \ot q_{ji}+1 \ot r_i.$$
   
   We note the following fact, which can be proved by adapting the 
 arguments of Section 5 of \cite{qiso_disc}, noting also that any compact interval can be transformed into the unit interval $[0,1]$ by dilation and translation.
          \blmma
      \label{interval_no_qiso}
      Let $I \subseteq \IR$ be a compact interval, $G$ the group of isometries of $I$ and $\gamma: C(I)\raro C(I \times G)
       $ be the canonical (co)action of $C(G)$. Let $\tau \in C(I)$ be given by $\tau(t)=t \forall t$, $\cls$ a unital 
       $C^*$ algebra, $\beta : C(I) \raro C(I) \ot \cls$ be a unital $\ast$-homomorphism such that $(\beta(\tau)(s)-\beta(\tau)(t))^2=
       (s-t)^2 1_\clq$ for all $s,t \in I$. Then there is a $\ast$-homomorphism $\pi: C(G) \raro \cls$ satisfying $\beta=({\rm id} \ot \pi)\circ \gamma$.
              \elmma
  
   \bthm
   \label{545}
   We  have $QISO(T_n,d) \cong \cls^+_n.$
   \ethm
    {\it Proof:}\\
      As $X_i$'s are self-adjoint, 
       $X_iX_j=0$ for $i \neq j$ and $\{ 1, X_1^2, \ldots, X^2_n\}$ are linearly independent, we get the following relations among the 
        $q_{ij}$ and $r_i$ :\\
        \be
        \label{1}
        r_i^*=r_i,~q_{ki}^*=q_{ki},~q_{ki}q_{kj}=0,~r_ir_j=0~~\forall k,~ i\neq j,\ee
        \be
        \label{2}
        r_iq_{kj}+q_{ki}r_j=0~\forall k~,i \neq j.\ee 
        Multiplying (\ref{2}) by the self-adjoint element $q_{ki}$  on the right we obtain \be \label{3} q_{ki}r_j q_{ki}=0.\ee
        But $r_i=\alpha(X_i)({\bf o})$  is a nonnegative element, hence (\ref{3}) implies $r_jq_{ki}=0=q_{ki}r_j$ for 
         all $i \neq j$ and all $k$. Moreover, using the isometry condition (\ref{pqr}) of Theorem \ref{iso_equiv} with $x,y \in T_j$ with $s,t$ at the $j$-th place respectively, 
         we have $\sum_{i=1}^n \left( X_j(s)-X_j(t) \right)^2  q_{ji}^2=(s-t)^2 1_\clq$ for all $s,t \in [0,1].$ This gives,
         \be \label{4} \sum_i q_{ji}^2=1 .\ee As $q_{ji}$ have mutually orthogonal ranges and are self-adjoint for fixed $j$, we get from the above that for fixed $j$, 
          $q^2_{ji}$ are mutually orthogonal 
          projections, say $\theta_{ji}$, with $\sum_i \theta_{ji}=1$.
          
          Next,
          for any fixed $i,j$, define $\beta_{ji} : C[0,1] \raro C[0,1] \ot \clq$ by $\beta_{ji}=\pi_j\circ \alpha_j \circ \pi_i^{-1},$ where $\alpha_j(f)=\alpha(f)_{T_j}$ 
           and $\pi_i : C^*(1,X_i) \raro C[0,1]$ is the isomorphism sending $X_i$ to the coordinate function (say $\tau$) of $[0,1]$.
           Thus, $\beta_{ji}(\tau)(t)=tq_{ji}+r_i$ for $t \in [0,1]$. This ( not necessarily coassociative)  unital  
            $\ast$-homomorphism satisfies $(\beta_{ji}(t)-\beta_{ji}(s))^2=(s-t)^2 1$ for all $s,t \in [0,1]$. It  follows from Lemma \ref{interval_no_qiso} that 
                          $r_i$ is also a projection and there is a projection $w_{ji} \leq \theta_{ji}=q^2_{ji}$ such that
             \be \label{5} q_{ji}=2w_{ji}-\theta_{ji},~~r_i=\theta_{ji}-w_{ji}.\ee
             Note that each $r_i$ is orthogonal to each $w_{jk}$ because for $k\neq i$, $w_{jk}$ is orthogonal to $\theta_{ji}$ which dominates $r_i$, and for $k=i$ the orthogonality 
              is contained in (\ref{5}). 
              
              We have $q_{ji}+q^2_{ji}=2 w_{ji}$, so  $\sigma(q_{ji}) \subseteq \{ t \in \IR:~t^2+t=0,2\}=\{ 0, 1, -1, -2\}$. But $-2$
               cannot be in the spectrum of $q_{ji}$ because $\sum_i q^2_{ji}=1$ implies $\| q^2_{ij}\| \leq 1$. Thus, $\sigma(q_{ji})$ is a subset of $\{0,1, -1\}$, 
               hence $q_{ji}(q^2_{ji}-1)=0$.
                Observe that $r_i=\frac{1}{2} (q^2_{ji}-q_{ji})$, so $r_i(q_{ji}+1)=\frac{1}{2}q_{ji}(q^2_{ji}-1)=0$. In other words, 
                $r_iq_{ji}=-r_i$, 
                 giving \be \label{zzz} (1 \ot r_i) \alpha(X_i)=-\sum_j X_j \ot r_i +1 \ot r_i.\ee Now, choose  $s,t \in (0,1)$  $j \neq k$ and let $x \in T_j$ and 
                  $y \in T_k$ be the elements having $s$ and $t$ as the $j$-th and $k$-th coordinates (elsewhere zero) respectively. 
                  The isometry condition gives, using also the orthogonality of different $r_i$'s and $r_i^2=r_i$ that $$ \left( r_i\alpha(X_i)(x)-r_i\alpha(X_i)(y) \right)^2=(s^2+t^2)r_i,$$ which, combined with (\ref{zzz}) and the fact that $r_i^2=r_i$ 
                   implies
                   $(s-t)^2 r_i=(s^2+t^2)r_i$, hence $r_i=0$. Thus, $w_{ji}=q_{ji}=\theta_{ji}$.

                  The  linear  span of $\{ X_1,\ldots, X_n \}$, say $\clw$, is left invariant by $\alpha$. 
                                         Using the arguments as in the proof of Theorem \ref{iso_equiv},  
                      we get an  $\alpha$-invariant state  $\phi$ (say) on $C(T_n)$ which is faithful on the unital $\ast$-algebra generated by 
                       the $X_i$'s and  the restriction of $\alpha$ to $\clw$ gives a unitary representation $U$ (say). 
                      However, $X_i=X_i^*$, $X_iX_j=0$ (hence $\phi(X_i^*X_j)=0$) for $i \neq j$ and the faithfulness of $\phi$ imply that 
                      $\{ c_i X_i, 1 \leq i \leq n \}$ gives a $\phi$-orthonormal basis, where $c_i>0, c_i^{-2}=\phi(X_i^2)$. 
                      The  matrix of $U$ w.r.t. this orthonormal basis 
                       is $(( \frac{c_i}{c_j} q_{ji} ))$. Thus, the antipode (say $\kappa$) gives: $\kappa(\frac{c_i}{c_j} q_{ji})
                       =\frac{c_j}{c_i}q_{ij}$ for all $i,j$ (as $q_{ji}^*=q_{ji}$). 
                       We already have 
                  $\sum_i q_{ji}= \sum_i \theta_{ji}=\sum_i w_{ji}=1.$ 
                       Applying $\kappa$ on this we obtain \be \label{7}
                        \sum_i \frac{c_j^2}{c_i^2} q_{ij}=1.\ee Moreover, for $i \neq k$,  $q_{ij}q_{kj}=\frac{c_i^2c_k^2}{c_j^4}\kappa(q_{jk}q_{ji})=0$.
                       Thus,  
                         $\{ q_{ij},i=1,\ldots,n\}$ is a family of mutually orthogonal  projections, hence (\ref{7}) implies 
                                                  $\frac{c_j^2}{c^2_i} q_{ij}=q_{ij}$ $\forall i$. But this  is possible only when either $q_{ij}$ is $0$ or $\frac{c^2_j}{c^2_i}=1$, 
                         which gives: \be \label{100} \sum_i q_{ij}=1.\ee 
                         
                         Thus, we get  a $C^*$ homomorphism from $\cls^+_n$ onto $\clq$ sending the canonical generators $w_{ij}$ to $q_{ij}$, which is clearly seen to be a CQG morphism. 
                           On the other hand, 
                            we get from \cite{huang1}  an action of $\cls_n^+$ on $C(T)$ which is clearly isometric in our sense,
                             hence (by the universality of $\clq$)  a surjective CQG morphism in the reverse direction. In other words, $\clq \cong \cls_n^+$.
                         \qed\\
              We next take up $T^\prime_n$. Recall from \cite{free_wr} the free wreath product of $\IZ_2$ by $\cls_n^+$, denoted by $\clu=\IC \IZ_2 \ast_w \cls_n^+$, 
               in particular the presentation  given on page 7 of that paper in terms of elements $\{ a_{ij} \}$. 
               Note that the isometry group of $[-1,1]$ is $\IZ_2=\IZ/2\IZ$ with the order $2$ generator of $\IZ_2$ giving the isometry $t \mapsto -t.$ 
               Using this as well as the natural action
                of $\cls_n^+$ on $C(T^\prime_n)$ as in \cite{huang1}, we can easily get an isometric action of $\clu$ which sends $X_i$ to $\sum_j X_j \ot a_{ji}$. Thus $\clu$ is a quantum subgroup 
                 of $\clq=QISO(T^\prime_n,d)$. We claim:
               \bthm
               \label{546}
           $QISO(T_n^\prime, d)$ is isomorphic with $\IC \IZ_2 \ast_w \cls_n^+$.  
               \ethm
    {\it Proof:}\\
          %
         Using the notation and arguments  of the proof of Theorem \ref{545}, we  can get $q_{ij}^*=q_{ij} \forall i,j$, 
         $q_{ki}q_{kj}=0$ for $i \neq j$ and for all $k$, and also $\sum_i q_{ji}^2=1.$
            Moreover, as in that proof, we consider the $\ast$-homomorphism $\beta_{ji}$ from $C[-1,1]$ to $C[-1,1] \ot \clq$ and conclude 
           by Lemma \ref{interval_no_qiso} that 
             $r_i=0$ and that $q_{ji}^2$ is a projection.  Thus, $\alpha(X_i)=\sum_j X_j \ot q_{ji}$ and using the mutual orthogonality of different $X_j$'s, 
            this gives
          $$ \alpha(X_i^2)=\sum_j X_j^2 \ot q^2_{ji}.$$ Using arguments of Theorem \ref{545} for proving  (\ref{100}) but replacing $X_i$ by $X_i^2$, 
          we can conclude $\sum_i q^2_{ij}=1$, as well as the orthogonality relations $q^2_{ik} q^2_{jk}=0$ for different $j,i$ and for all $k$. 
          But $q_{ij}^2=q_{ij}^* q_{ij}$ is a projection, so that $q_{ji}$ is a normal partial isometry with the domain and range projections being $q^2_{ji}.$ 
          Therefore, orthogonality of $q^2_{ik}$ and $q^2_{jk}$ for different $i,j$ imply the orthogonality of $q_{ik}$ and $q_{jk}$ too. 
          This completes the proof that $q_{ij}$'s satisfy the same relations as the generators $a_{ij}$'s of  $\clu$ given in \cite{bichon},
          which gives a morphism from $\clu$ to $\clq$ sending $a_{ij}$ to $q_{ij}$ and thereby completing the proof. \qed\\
    \brmrk
    $QISO(T_n,d)$ for $n\geq 4$ and $QISO(T_n^\prime,d)$ for all $n \geq 2$ are genuine CQG's. 
        \ermrk
        \brmrk
        Let $(X_n, d_n)$ be the topological join at a common point $x_0$ of 
      $n$ copies of a compact metric space $(X,d)$ isometrically embedded in some $\IR^m$ and let $\clq_{x_0}$ be
       the maximal quantum subgroup of $QISO(X,d)$ which fixes the point $x_0$ in some suitable sense. 
      Theorem \ref{545} and Theorem \ref{546} suggest the following generalization of  results in \cite{bichon_banica}:
      $QISO(X_n, d_n) \cong \clq_{x_0} \ast_w \cls_n^+ $.
             \ermrk
 
{\bf Acknowledgment:}\\
I thank an anonymous referee for valuable suggestions and comments leading to substantial improvement of the exposition, 
 including part of the discussion in the Introduction 
 about motivations of formulation of quantum isometry groups for metric spaces. I
 gratefully acknowledge  A. L. Chirvasitu for several useful discussions through email, in particular about the issue of injectivity of actions and the fact that 
  our main result applies to all finite metric spaces. Several other mathematicians also deserve special mention and grateful acknowledgement:
   A Skalski, S. Wang, J. Bhowmick, H. Huang, M. Rieffel, J. Quaegebeur, M. Sabbe and P. Etingof for their numerous comments, suggestions, corrections etc., 
   R. S. Hajra for invaluable help in drawing figures and S. Vaes and M Rieffel for invitations to KU Leuven and UC Berkeley (respectively).

\end{document}